\documentclass{amsart}
\usepackage{amsfonts}

\newcommand{\C}{{\mathcal C}}
\newcommand{\M}{{\mathcal M}}
\newcommand{\R}{{\Bbb R}}
\newcommand{\N}{{\Bbb N}}
\theoremstyle{plain}
\newtheorem{theorem}{Theorem}
\newtheorem{corollary}{Corollary}

\newtheorem{lemma}{Lemma}

\theoremstyle{definition}

\newtheorem{remark}{Remark}
\newtheorem{example}{Example}
\begin{document}

\title[global stability criterion]{Yorke and Wright $\mathbf{3/2}$-stability
theorems  from a unified point of view}

\thanks{This research was supported by FONDECYT (Chile), project 8990013.
E. Liz was supported in part by M. C. T. (Spain) and FEDER under project BFM2001-3884.
V. Tkachenko was supported in part by F.F.D. of Ukraine, project 01.07/00109.}

\author[E. Liz,  V. Tkachenko, S. Trofimchuk]{}
\email{eliz@dma.uvigo.es} 
\email{vitk@imath.kiev.ua}
\email{trofimch@uchile.cl}


\date{\today}

\maketitle

\centerline{\scshape  Eduardo Liz}
\medskip

{\footnotesize
\centerline{ Departamento de Matem\'atica Aplicada II, E.T.S.I. Telecomunicaci\'on}
\centerline{ Universidad de Vigo, Campus Marcosende, 36280 Vigo, Spain}
}

\medskip
\centerline{\scshape  Victor Tkachenko}
\medskip

{\footnotesize
\centerline{ Institute of Mathematics,
National Academy of Sciences of Ukraine}
\centerline{ Tereshchenkivs'ka str. 3, Kiev, Ukraine}
}

\medskip
\centerline{\scshape  Sergei Trofimchuk}
\medskip

{\footnotesize
\centerline{ Departamento de Matem\'aticas, Facultad de Ciencias }
\centerline{ Universidad de Chile, Casilla 653, Santiago, Chile}
}

\bigskip
\begin{quote}{\normalfont\fontsize{8}{10}\selectfont
{\bfseries Abstract.}
We consider a family of scalar delay differential equations
$x'(t)=f(t,x_t)$, with a nonlinearity $f$
satisfying a  negative feedback condition combined with
a boundedness condition.
We present a global stability criterion for this family, which in particular unifies the
celebrated
$3/2$-conditions given for the Yorke  and the Wright type
equations. We illustrate our results with some applications.
\par}
\end{quote}

\section{Introduction and main result}
\label{sec.int}

In this paper we present a global stability criterion for a family of scalar
functional differential equations
\begin {equation} \label{f}
x'(t) = f(t,x_t), \ \  (x_t(s) \stackrel{def}=
x(t+s), \ s \in [-1,0]),
\end{equation}
where $f: \Bbb R \times \C \to \Bbb R$ is a measurable
functional, $\C=\C([-1,0], \R)$

As we will show in Section 3, our setting allows us to prove global stability
results for a large family of functional delay differential equations, including a
very general form of the delayed logistic equation, and differential equations with
maxima among others.

Next we introduce the hypotheses that  will be required in Eq. (\ref{f}).  In
order to understand the motivation for the choice of these conditions, we recall
some classical results. We refer to the famous $3/2$
stability results due to Myshkis \cite{my}, Wright \cite{wr} and Yorke
\cite{yo}.

In particular, the Wright equation can be written in the
form
\begin {equation} \label{w}
x'(t) =   f(x(t-1)),
\end{equation}
with $f(x)=p(e^{-x}-1)$, $p>0$. The famous  $3/2$ stability result by Wright says
that all solutions of this equation converge to zero if $p\leq 3/2$. The closeness
between this condition and the (local) asymptotic stability condition, $p<\pi/2$,
suggests the equivalence between the local and the global asymptotic stability
(this is the famous Wright's conjecture, which still remains open).

However, it is known   that number $3/2$ is the best bound when
we consider differential equations with variable delay , even in the linear case
\begin {equation} \label{l}
x'(t) =  -px(t-h(t))\, , \, p>0,\, 0<h(t)\leq h.
\end{equation}
Myshkis proved that Eq. (\ref{l}) is exponentially stable if $ph<3/2$, and it is
possible to find examples such that $ph=3/2$ and (\ref{l})
has a nontrivial periodic solution (see \cite[p.~191]{yo}).

In 1970, Yorke extended the Myshkis criterion to a family of scalar functional
differential equations (\ref{f})
where $f:\R\times\C([-h,0])\to\R,\, h>0,$  is continuous and satisfies the following conditions:
\begin{enumerate}
  \item[{\rm \bf(Y1)}] There exists $a < 0$ such that
$$ a\mathcal{M}(\phi) \leq f(t,\phi)\leq
-a\mathcal{M}(-\phi)
$$
for all $\phi \in \C$, where
$\mathcal{M}(\phi) = \max\{0, \max_{s\in[-h,0]}\phi(s)\}.$
  \item[{\rm \bf(Y2)}]  For all sequences $t_n\to \infty$ and $\phi_n$
converging to a constant nonzero function in $C$, $f(t_n,\phi_n)$ does not
converge to $0$.
\end{enumerate}

Under these conditions, if $0<|a|h<3/2$ then all solutions of (\ref{f}) converge
to zero as $t\to \infty$.

We notice that condition {\rm \bf(Y2)} 
 is only required to guarantee that the solutions of (\ref{f}) that
monotonically converge to a constant in fact should converge to zero.

One can check that condition {\rm \bf(Y1)} is not satisfied by the Wright
equation, and therefore Wright's theorem cannot be deduced from the Yorke result.

Trying to generalize the Wright theorem, in \cite{lprtt} we prove the following
result:
\begin{theorem}
\label{t1}
Assume that $f\in \C^3(\R,\R)$ and it satisfies the following conditions:
\begin{enumerate}
 \item[{\bf (A1)}] $xf(x) < 0$  for  $x\neq 0$ and $f'(0) < 0$.
 \item[\bf (A2)] $f$ is bounded below and has at most one critical point
$x^*\in \R$ which is a local extremum.
 \item[\bf (A3)] $(Sf)(x)<0$ for all $x\neq x^*$, where
$
Sf=f'''(f')^{-1}-(3/2)(f'')^2(f')^{-2}
$
is  the  Schwarz derivative of $f$. 
\end{enumerate}
If  $|f'(0)| \leq 3/2$,
then the steady state
solution $x(t) = 0$ of Eq. (\ref{w}) is globally
attracting.
\end{theorem}

\begin{remark} Conditions {\bf (A1)}-{\bf (A3)} are satisfied by the Wright
equation and other more complicated equations arising in population dynamics (see
\cite{lprtt} for details).
\end{remark}

Conditions {\bf (A1)}-{\bf (A2)} are not sufficient for the global attractivity in
(\ref{w}) (see \cite{lprtt, waltdcds}).
Hence an additional condition is required. We note that  condition {\bf (A3)} is not
the unique option, 
in fact we  only need  some geometric consequences of the inequality $Sf < 0$
for the graph of $f$. In particular,  the following key
result (\cite[Lemma 2.1]{lprtt}) is very important:
\begin{lemma}
\label{l1}
Assume that $f$ satisfies {\bf (A1)}-{\bf (A3)}, and $f''(0)>0$. Let $a=f'(0)<0$,
$b=-f''(0)/(2f'(0))>0$, and $r(x)=ax/(1+bx)$. Then $r(x)>f(x)$ if $x\in (-1/b,0)$, and
$r(x)<f(x)$ if $x>0$.
\end{lemma}

From Lemma \ref{l1}, we can obtain immediately the following

\begin{corollary}
\label{c1}
Assume that $f$ satisfies {\bf (A1)}-{\bf (A3)}, and let $f(t,\phi)=f(\phi(-1))$.
Then the following ``{\it generalized Yorke condition}":
\begin{enumerate}
  \item[{\rm \bf(GY)}] There exist $a< 0, \ b\geq 0$ such that
$$ r(\mathcal{M}(\phi))=\frac{a\mathcal{M}(\phi)}{1+b\mathcal{M}(\phi)} \leq
f(t,\phi)\leq
\frac{-a\mathcal{M}(-\phi)}{1-b\mathcal{M}(-\phi)}=r(-\mathcal{M}(-\phi)),
$$
\end{enumerate}
where the first inequality holds for all $\phi\in\C$, and the
second one for all
$\phi$ such that
$\min_{s\in[-1,0]}\phi(s)> -b^{-1}
\in [-\infty,0)$.
 Here
$\mathcal{M}(\phi) = \max\{0, \max_{s\in[-1,0]}\phi(s)\}$ is the Yorke functional.
\end{corollary}

This corollary suggests the unification of Yorke's and Wright's $3/2$ results by
using condition {\rm \bf(GY)}. In fact, we introduce the following three hypotheses
{\rm \bf(H)}:
\begin{enumerate}
  \item[{\rm \bf(H1)}] $f: \Bbb R \times \C \to \Bbb R$ satisfies the Carath\'eodory condition.
Moreover, for every  $q \in \Bbb R$  there exists  \ $\vartheta(q)
\geq 0$ such that $f(t,\phi) \leq \vartheta(q)$  almost
everywhere on $\Bbb R$ for every $\phi \in \C$ which satisfies the
inequality $\phi(s) \geq q, s \in [-1,0]$.
  \item[{\rm \bf(H2)}]  Condition {\rm \bf(GY)} holds.
  \item[{\rm \bf(H3)}]  $\int_0^{+\infty}f(s,p_s)ds$ diverges for every continuous
$p(s)$ having nonzero limit at infinity.
\end{enumerate}
We recall that $f(t,\phi)$ is a Carath\'eodory function if it is measurable in $t$ for each fixed
$\phi$, continuous in $\phi$ for each fixed $t$, and for any fixed $(t,\phi)\in\R\times\C$ there
is a neighbourhood $V(t,\phi)$ and a Lebesgue integrable function $m$ such that
$|f(s,\psi)|\leq m(s)$ for all $(s,\psi)\in V(t,\phi)$ (see
\cite[p.~58]{hl}). 

The following result
improves the above mentioned theorems due to Wright and Yorke
respectively (notice that {\rm \bf(H3)} implies that $x(t) \equiv 0$ is
the unique equilibrium of the equation):

\begin{theorem}
\label{main}
Assume that $f$ satisfies {\rm \bf(H)} and either $b>0$ and $|a| \leq 3/2$, or
$b=0$ and $|a| < 3/2$. Then all  solutions of (\ref{f}) converge to zero as
$t\to +\infty$.
\end{theorem}

\begin{remark}
\begin{enumerate}
\item If {\rm
\bf (H2)} holds with $b=0$ ({\it Yorke condition}), then {\rm \bf (H1)} is satisfied
automatically with $\vartheta(q) = -a\M(-q)$.
\item Conditions {\rm \bf (H1)}- {\rm \bf (H3)} are  satisfied for Eq. (\ref{w})
under our hypotheses {\rm \bf(A1)}-{\rm \bf (A3)}, with $a=f'(0)$,
$b=-f''(0)/(2f'(0))$.
\item 
The constant $3/2$ in Theorem \ref{main} is the best possible. For $b=0$, this
sharpness was shown in \cite{yo} (see also \cite{ilt, ltt}). For $b>0$, Theorem \ref{main} can be
applied to prove that all positive solutions of the logistic equation with variable delays
\begin{equation}
\label{levd}
x'(t)=px(t)(1-x(t-h(t)))\, ,\, 0<h(t)\leq h,
\end{equation}
converge to $1$ if $ph\leq 3/2$ (see Theorem \ref{apl} in Section 3). Since the
linearized equation of (\ref{levd}) at $x=1$ is (\ref{l}), constant $3/2$ cannot be
improved.  
It is a remarkable fact that when  Yorke's result cannot be
extended to the value $a=-3/2$ for $b=0$, Theorem \ref{main} allows this for every
$b>0$.
\end{enumerate}
\end{remark}

If $b>0$, since $\mathcal{M}$ is a positively homogeneous functional
($\mathcal{M}(k\phi) = k\mathcal{M}(\phi)$ for every $k\geq 0$,
$\phi \in \C$), and since the global attractivity property of the
trivial solution of (\ref{f}) is preserved under the simple
scaling $x = b^{-1}y $,  the exact value of $b$ does not have
importance and we can assume that $b = 1$. Also, the change of
variables $x = -y$ transforms (\ref{f}) into $y'(t) =
-f(t,-y_t)$ so that it suffices that at least one of the two
functionals $f(t,\phi), -f(t,-\phi)$ satisfies {\rm \bf(GY)}.

\section{Proof of Theorem \ref{main}.}
\subsection{Auxiliary results.} Throughout this subsection, we will assume that
$b=1$ (and hence $r(x)=ax/(1+x)$).

\begin{lemma} \label{L1} Let {\rm \bf(H)} hold and
$x:[\alpha -1,\omega) \to \R$ be a solution of (\ref{f})
defined on the maximal interval of  existence. Then $\omega= +
\infty$ and $M = \limsup_{t \to \infty}x(t)$, $ m = \liminf_{t \to
\infty}x(t)$ are finite. Moreover, if $m \geq 0$ or $M \leq 0$,
then $M = m  =0$.
\end{lemma}
\begin{proof} Note that {\rm \bf(GY)} implies that $f(t,\phi) \geq a$ for all $t \in \R$
and $\phi \in \C$.
We claim that every solution $x(t,\gamma)$ with  initial value $\gamma$ such that
$ q \leq \gamma(s) \leq Q$, $s\in [-1,0],$  satisfies the inequality
$$
x(t,\gamma) \geq \min \{q,0\} + a= \kappa, \quad t \geq 0.
$$
Indeed, if there is $\delta > 0$ such that
$x(t,\gamma)= \min \{q,0\} + (1+\delta)a $ for the first time at some point  $t=
u\geq 0$, then, for every $w \in (u-1-\delta/2,u)$,
$$
x(u,\gamma)-x(w,\gamma)= \int_{w}^{u}f(s,x_s(\gamma))ds \geq a(u-w) >
a(1+\delta/2).
$$
Hence $x(w,\gamma)\leq a\delta/2 <0 $  for all $w \in (u-1-\delta/2,u)$, and
therefore  $\M(x_s)=0$ for all $s\in (u-\delta/2,u)$. Thus, by {\rm \bf(H2)}, $x'(s,\gamma)=
f(s, x_s(\gamma))\geq 0$ within some left neighborhood of $u$, contradicting  the choice of
this point.

Proceeding analogously and using {\rm \bf(H1)}, we obtain that
$$
x(t,\gamma) \leq \max \{Q,0\}+ \vartheta (\kappa), \quad t \geq 0.
$$
Hence $x(t)$ is bounded on the maximal interval of existence that
implies the boundedness of the right hand side  of Eq. (\ref{f}) along
$x(t)$.  Thus $\omega = +\infty$ due to the corresponding continuation
theorem (see \cite{hl}).

Now, suppose that $M < 0$. Then $x(t) < M/2 < 0$ beginning from some
$t' = d-1$ so that, in view
of $x'(t) \geq 0$, the solution $x(t)$ converges monotonically to the
negative value $x(+\infty) = M$. We get a contradiction since, by {\bf (H3)},
$$ x(t) = x(d) + \int_d^t f(s,x_s)ds
 \to +\infty.$$
 Next we   consider  the situation when
$M=0$ and  $m < 0$. In this case, $x(t)$ necessarily oscillates  about
zero. Indeed, otherwise $x(t) \leq 0$ and thus $x'(t) = f(t,x_t) \geq 0$,
so that $x(t)$ converges monotonically to the trivial steady
state (implying $m =0$).  Now, since $x(t)$ is oscillating, we
can find a sequence of intervals $I_k= (l_k,r_k)$ containing $e_k$
such that $x(t) < 0,\ t \in I_k$ and  $\min_{I_k} x(t) = x(e_k) \to m$ as $ k \to
+\infty$, while $e_k$ is the minimal point from $I_k$ having this property.
We claim that $e_k - l_k \leq 1$. On the contrary, let us suppose
that $e_k - l_k > 1$. Then $x_t < 0$ (and, consequently, $x'(t) \geq 0$)
for all $t$ from a small neighborhood of $e_k$, contradicting to the
choice of $e_k$ as the  leftmost  point of global minimum in $I_k$. Finally,
observing that
$$x(e_k) = \int_{l_k}^{e_k}f(s,x_s)ds \geq  r(\max_{u \in [l_k-1,l_k]} x(u)) \to 0,
\ k \to +\infty,$$
we get a contradiction with the relation $x(e_k) \to m < 0$.

The case $m\geq 0$ is similarly addressed.
\end{proof}
\begin{remark}\label{rema} The last part of the above proof can be repeated to
analyze the structure of the set of extreme points for every solution
$x(t)$ satisfying $m < 0$ and $M >0$. We see that in that case there
exist  sequences of  intervals $A_k = (a_k,b_k), \ A'_k = (a'_k,b'_k)$
and points $e_k \in A_k, \ e'_k \in A'_k$ such that $x(a'_k) = x(a_k) =0, \
e_k-a_k \leq 1, e'_k -a'_k \leq ,1$ while $x(e_k) \to m, \ x(e'_k) \to  M$
and $x(t)$ does not change  sign over $A_k, A'_k$. Moreover, for
each $k$, $e_k$ and $e'_k$  could be chosen as the points of global maximum of
$|x(t)|$ on  $A_k$ and $A'_k$ respectively. 
\end{remark}

Now, we define  continuous functions $A,B:
(-1,+\infty)\to \R$ and $D:\R_+ \to \R$ by
$$
A(x)=x+r(x)+\frac{1}{r(x)}\int_{x}^{0}r(t)dt,
\quad
B(x)=\frac{1}{r(x)}\int_{-r(x)}^{0}r(s)ds
\quad {\rm for \ } x\not=0,
$$
$$
A(0) = B(0) = 0, \quad
D(x)= \left\{ \begin{array}{ll} A(x) & \textrm{if $r(x)<-x$},
\\ B(x) & \textrm{if $r(x)\geq -x$}. \end{array} \right.
$$
For the case $a <-1$, we will also use  the rational function
$$R(x) = \frac {(A'(0))^2\ x}{A'(0) - (A''(0)/2)x}$$
defined on the interval
$(2A'(0)/A''(0), \infty) = (\nu,\infty)$.  Note that $A'(0)= a +1/2<0,
\ A''(0) = -(2a + 1/3) >0$. It is easy
to check that $a < -1/6$ implies $(\nu,+\infty) \subset (-1, +\infty)$,
and that $A(x_2)=B(x_2)$, where $r(x_2) = -x_2 <0$. Also $B'(0)= -(r'(0))^2/2 =
-a^2/2$.

The following relations were established in \cite[Lemma 2.4, Corollary 2.7]{lprtt}:

\begin{lemma}
\label{32}
If $a < -1,$ then
$(A(x)-R(x))x > 0$
for  $x \in (\nu,x_2) \setminus \{0\}$.
Moreover, if $a   \in [-1.5, -1.25]$, then $D(x) > R(x)$ for $x > 0$.
\end{lemma}

In the next stage of the proof, we establish various relations
between $m$ and $M$, all of them being expressed in terms of the recently
introduced functions. Notice that, by Lemma \ref{L1}, the only
case of interest is when $m<0<M$; thus we can suppose the existence
of points $t_j, s_j$ of local maxima  and local minima respectively
such that $x(t_j)= M_j \to M, x(s_j) = m_j \to m$ and $s_j, t_j \to
+\infty$ as $j \to \infty$.

\begin{lemma}
\label{44}
We have $m \geq D(M)$ and $m \geq r(-r(M)/2)$.
\end{lemma}
\begin{proof}
First we assume that $r(M) < -M$. Then
$M/r(M) \in (-1,0]$.
Take now $\varepsilon > 0$ such that
$t_1 = (M + \varepsilon)/r(M + \varepsilon) \in (-1,0]$.
Obviously, $m_j > m-\varepsilon$ and $M_j < M+\varepsilon$ for
sufficiently large $j$.  By Remark \ref{rema},
there exists $\tilde s_j \in [s_j - 1, s_j]$
such that $x(\tilde s_j) = 0$ and $x(t) < 0$ for
$t \in (\tilde s_j, s_j].$

Next,  $z(t)=r(M+\varepsilon)(t - \tilde s_j), \
t \in [\tilde s_j + t_{1}, \tilde s_j + t_{1}+1]$
solves the initial value problem $z(s)=M+\varepsilon, \
s \in [\tilde s_j + t_1-1, \tilde s_j + t_1]$
for
\begin{equation}
\label{ep}
z'(t)=r(z(t-1)).
\end{equation}

Clearly $M+\varepsilon=z(t) > x(t)$ for all
$t \in [t_{1} + \tilde s_j - 1, \tilde s_j + t_{1}]$.
Moreover, we will prove that $z(t) > x(t)$ for all
$t \in [\tilde s_j + t_{1}, \tilde s_j)$.

Indeed, if this is not the case we can find
$t_{*} \in [\tilde s_j + t_{1}, \tilde s_j)$ such that
$z(t_{*})=x(t_{*})$ and
$z(t)> x(t)$ for all $t \in [\tilde s_j + t_1, t_{*}]$.
We claim  that
\begin{equation}
\label{zvx}
x'(t)> z'(t)\ \mbox{ for all }\  t \in [t_{*}, \tilde s_j].
\end{equation}
We have
$$
z'(t)=r(z(t-1)) < r(\mathcal{M} (x_t)) \leq f(t, x_t) =x'(t).
$$
 After integration
over $(t_*, \tilde s_j)$, it follows from (\ref{zvx}) that
$x(t_*)<z(t_*)$, which is a contradiction.

Thus $z(t) > x(t)$ for $t \in [t_{1} + \tilde s_j, \tilde s_j)$
and, arguing as above, we obtain
 \begin{eqnarray*}
m_j &=& \int_{\tilde s_j}^{s_j} x'(s)ds =
\int_{\tilde s_j}^{s_j} f(s, x_s) ds \geq
\int_{\tilde s_j}^{s_j} r(\mathcal{M} (x_s)) ds > 
 \int_{\tilde s_j}^{s_j} r(z(s-1)) ds\\ &=&
\int_{\tilde s_j - 1}^{\tilde s_j + t_1} r(M + \varepsilon) ds +
\int_{\tilde s_j + t_1}^{s_j - 1} r(r(M + \varepsilon)(s - \tilde s_j))ds \\
&= &  r(M + \varepsilon)(t_1 + 1) +
\int_{t_1}^{s_j - \tilde s_j - 1} r(r(M + \varepsilon)u)du  \\
& \ge&  M + \varepsilon + r(M + \varepsilon)+  \int_{t_{1}}^{0}r(r(M+\varepsilon)u)du
= A(M+\varepsilon).
\end{eqnarray*}

As a limit form of this inequality, we get $m\geq A(M)$.

In the  general case (i.e. we do not assume that $r(M) < -M$), we use
the inequality $f(t, x_t) \ge r(\mathcal M (x_t)) > r(M+\varepsilon)$ to
see that, for $t \in (\tilde s_j-1, \tilde s_j)$,
$$
x(t) = -\int_t^{\tilde s_j} x'(s)ds < -\int_t^{\tilde s_j} r(M+\varepsilon)ds =
r(M+\varepsilon)(t - \tilde s_j).
$$
In consequence,
\begin{eqnarray*}
m_j&=&x(s_j) =\int_{\tilde s_j}^{s_j}f(s, x_s)ds >
\int_{\tilde s_j}^{s_j}r(r(M+\varepsilon)(s- \tilde s_j - 1))ds  
\\
&\ge& \int_{0}^{s_j - \tilde s_j}r(r(M+\varepsilon)(s-1))ds \ge
\int_{0}^{1} r(r(M+\varepsilon)(u-1)) du = B(M+\varepsilon).
\end{eqnarray*}
Therefore, we obtain that $m\geq B(M)$.
Finally, applying Jensen's integral inequality (see \cite[p. 110]{roy}) , we have
$$m \geq B(M) =\frac{1}{r(M)}\int_{-r(M)}^{0}r(s)ds \geq  r(-r(M)/2). $$
This completes the proof.
\end{proof}
As a consequence of Lemmas \ref{32}, \ref{44}, we obtain that
$R(m),r(m)$ and $r(r(-r(M)/2))$ are well-defined and that
$R(\nu,+\infty) \subset (\nu,+\infty)$ for  suitable values of
$a$:
\begin{corollary}
\label{cy}
\hspace{-2mm} We have $m > -1, r(-r(M)/2) >-1$ if
$a \in [-1.5,0)$ and $m > \nu$ if $a\in [-1.5,-1.25]$.
\end{corollary}
\begin{proof} Indeed, for  $a \in (-2,0)$ we have
$$
m \geq r(-r(M)/2) > r(-\frac{r(+\infty)}{2})= \frac{-a^2}{2-a} \geq -1.
$$
Next, $-A'(0) = -(a+0.5) \leq 1$ for $a \geq -1.5$, and
Lemmas \ref{32}, \ref{44} lead to the estimate
$$m  > D(+\infty) = B(+\infty) \geq R(+\infty) = -A'(0)\nu \geq \nu.$$
This proves the corollary.
\end{proof}
\begin{lemma}
\label{43} Let $a \in [-1.5,0)$. We have
$M \leq r(m)$.  Moreover, if $a \in [-1.5,-1.25]$ then
$M <  R(m)$.
\end{lemma}
\begin{proof} We have that $r(m)$ is well defined and $[m, +\infty) \subset [-1,
+\infty)$ since $a \in [-1.5,0)$ (see Corollary \ref{cy}).
Let us consider a solution $x(t)$ of Eq. (1)
and take $s_j, t_j, m_j, M_j,$ as in the  paragraph below
 Lemma \ref{32}. By Remark \ref{rema}, there exists $\tilde t_j \in [t_j - 1, t_j)$
such that $x(\tilde t_j) = 0.$

First we note that, for $\varepsilon > 0$ and $j$ sufficiently large, $f(s, x_s)
\le r(-{\mathcal M} (-x_s)) < r(m-\varepsilon)$  for  $s\in [\tilde t_j - 1, \tilde
t_j]$
 and $a \in [-1.5,0)$. Thus
$$
M_j = x(t_j) = \int_{\tilde t_j}^{t_j}f(s, x_s)ds <
\int_0^1r(m-\varepsilon)ds = r(m-\varepsilon).
$$
By taking the limits as $\varepsilon\to 0$ and $j\to+\infty$, we obtain the inequality
$M\leq r(m)$.

 Now, if $a<-1$ then  $r(m)>-m$, from which for all sufficiently small
$\varepsilon >0$ we obtain that
$t_2=(m- \varepsilon) (r(m-\varepsilon))^{-1} \in (-1,0]$.  Next,
$z(t)= r(m-\varepsilon)(t - \tilde t_j)$, with $\ t \in [\tilde t_j + t_2, \tilde t_j + t_2 +1]$,
solves the initial value
problem $z(s)=m-\varepsilon,\ s \in [\tilde t_j + t_2 - 1, \tilde t_j + t_2]$
for Eq. (\ref{ep}). Now we only have to argue as in the proof of Lemma \ref{44}
to find out that
$$M_j = \int_{\tilde t_j}^{t_j} f(t, x_t) dt \leq
\int_{\tilde t_j}^{t_j} r(-\mathcal{M}(-x_t)) dt \le
\int_{\tilde t_j}^{t_j}r(z(t-1)) dt \leq A(m-\varepsilon). $$

Thus $M\leq A(m)$. Finally, by Lemma \ref{32} and Corollary \ref{cy}, we obtain
$M\leq A(m)<R(m)$ when $a\in [-1.5,-1.25]$.
\end{proof}

\subsection{Proof of Theorem \ref{main}}
 Let $x:[\alpha-h, +\infty) \to \R$ be a solution of Eq.
(\ref{f}) and set $M = \limsup_{t \to \infty}x(t), \ m =
\liminf_{t \to \infty}x(t)$. We will reach a contradiction if we
assume that $m<0<M$ (note that the cases $M \leq 0$ and $m \geq 0$
were already considered in Lemma \ref{L1}).

Assume first that $b=1$.
If $a \in (-1.5,0)$, in view of Lemmas \ref{44}, \ref{43} and Corollary \ref{cy}
we obtain that $M \leq r(m) \leq r\circ r(-r(M)/2) = \lambda (M)$ with the
rational  function $y =\lambda(x)$. Now, $\lambda(M) < M$ for $M > 0$
if $\lambda'(0) =(1/2) |a|^3< 1$. Therefore,  if $a \in [-1.25, 0)$ we obtain
the desired contradiction under the assumption $M>0$.

Now let $a \in [-1.5, -1.25]$ and, consequently,
$R'(0) = a + 0.5 \in [-1,-0.75]$.
In this case Lemmas \ref{32}, \ref{44} and \ref{43} imply that $M < R(R(M))$.
As $R\circ R(x) \leq  x$ for all $x > 0$ whenever
$(R\circ R)'(0) = (R'(0))^2 \leq 1$, we obtain a contradiction again.
 Therefore $x(t)\equiv 0$ is the global attractor of
Eq. (\ref{f}) if $a \in [-1.5,0)$.

If $b=0$, we employ the linear function $r(x)=ax$. Arguing as in the proofs of
Lemmas \ref{44} and \ref{43}, we obtain that $m\geq (a+1/2) M$, $M\leq (a+1/2) m$
for
$a\in (-3/2, -1]$, and $m\geq (-1/2)a^2 M$, $M\leq am$ for $a<0$ (see \cite{lprtt}
for more details).

Hence, if $a\in (-3/2, -1]$, we get the contradiction $M\leq (a+1/2)^2 M<M$. Finally, if $a\in (-1,0)$, we obtain $m\geq
(-1/2)a^2 M\geq (-1/2) a^3 m$. Thus, $-a^3/2\geq 1$, a contradiction.

\section{An application}

Probably, the most interesting object to which we can apply our
results is the following generalization of the logistic delayed
equation:
\begin{equation}
\label{gle} x'(t) = \lambda(t)x(t)f(t,\mathcal{L}(t,x_t)),\quad \
\quad t\geq 0,\ x \geq 0.
\end{equation}
Here $\lambda: [-h, \infty)  \to 
(0, \infty)$ is measurable and 
$\int^{\infty}_0 \lambda (s) ds = \infty, 
\sup_{t \geq 0}\int_{t-h}^{t} \lambda (s)ds < \infty.$ 
We suppose that $f: \R_+ \times \R \to \R $
and $\mathcal{L}(t, \phi): \R_+ \times \C([-h,0]) \to \R$ are 
Carath\'eodory functions, and that 
$$\min \phi \leq \mathcal{L}(t, \phi) \leq \max \phi$$ 
for every  $\phi \in \C([-h,0]), \ t \geq 0$. 
We are interested in the case when, apart from
$x(t) \equiv 0$, Eq. (\ref{gle}) has another equilibrium $x(t)
\equiv \kappa$. Without  loss of generality we can assume that 
$\kappa =1$ (and, consequently, that $f(t,1) \equiv 0$). Finally, we
will require the divergence of  $\int_0^{+\infty}
\lambda(s)f(s,w(s))ds$ for every continuous function $w$ converging to
some positive number different from $1$, as well as the
following negative feed-back condition:
$$ x\ f(t,1+x)< 0\, ,\forall\,x >
-1, \ x \not= 0.
$$
 As a simple observation shows, 
every nontrivial solution of (\ref{gle}) is eventually positive, so that 
we will only study the behavior of the positive 
solutions.

The next result is a consequence of Theorem \ref{main}:
\begin{theorem}
\label{apl}  Assume that 
\begin{equation}
\label{nf}
x\ (f(t,x + 1)-r(x)) \geq 0
\end{equation}
for some $r(x) = ax/(1+bx)$ with $ a < 0,\ b \geq 0$ and for
all $x
> \max\{-1,- b^{-1}\}$. 

If  $b \not= 0.5$
and, for some $T > 0$, $|a| \Lambda  \leq
3/2$,  
  then $\lim_{t \to +\infty} x(t) = 1$ for every 
nontrivial solution of (\ref{gle}), where $\Lambda = \sup_{t\geq T} \int_{t-h}^{t}
\lambda(s)ds$.  If $b =0.5$, then the same conclusion holds when  \ $|a|\Lambda <
3/2$.
\end{theorem}
\begin{proof} 
First, let $b \geq 0.5$. The change of
variables $y(s) = \ln x(t)$, where $s = s(t) = \Lambda
^{-1}\int_0^t\lambda(u)du$ (with the inverse $t= t(s)$ to
$s=s(t)$), reduces (\ref{gle}) to
\begin{equation}
\label{glem} y'(s) = \Lambda g(s,\mathcal{K}(s,y_s)),\quad \ y \in
\R, \quad s\geq 0.
\end{equation}
Here, $g(s,x) = f(t(s),x)$  and 
$\mathcal{K}(s,\phi) = \mathcal{L}(t(s),\Psi)$, where
$\Psi=\Psi(s,\phi)$ is defined by $[\Psi(s,\phi)](u)=\exp[
\phi(- \sigma(s,u))]$, being $\sigma(s,u) = \Lambda^{-1}
\int_{t(s)+u}^{t(s)}\lambda(v)dv$.  Notice that $\mathcal{K}: \R_+ \times \C([-1,0])
\to \R$ is a Carath\'eodory  function, $\exp(-\mathcal{M}(-\phi)) 
\leq \mathcal{K}(s, \phi) \leq \exp(\mathcal{M}(\phi))$ 
for every  $\phi \in \C([-1,0]), \ s \geq 0$, and
$\mathcal{K}(s,0) \equiv 1$. Since
$g(s,\mathcal{K}(s,\phi)) = f[t(s),1 + (\mathcal{ K}
(s, \phi) -1)]$ we have, for $\mathcal{K} (s,
\phi) \leq 1$,  that 
\begin{eqnarray}\label{ar1}
r_1(\mathcal{M}(\phi)) &\leq& 0 \leq  \Lambda
g(s,\mathcal{K}(s,\phi)) \leq \Lambda r(\mathcal{K}
(s,\phi)-1)\\ &\leq&   \Lambda r(\exp
[-\mathcal{M}(-\phi)]-1)  \leq r_1(-\mathcal{M}(-\phi)), \nonumber
\end{eqnarray}
where $r_1(x) = a\Lambda x/(1+(b-0.5)x).$ Indeed, function $v(x)=r(e^x-1)$ satisfies
conditions {\bf (A1)}-{\bf (A3)}
in Theorem \ref{t1} and therefore, by Lemma \ref{l1}, $r(e^x-1)\leq r_1(x)$ for
$x<0$, and $r(e^x-1)\geq r_1(x)$ for
$x>0$, since $v'(0)=a$, $v''(0)=a(1-2b)$.

Analogously, if
$\mathcal{K} (s,\phi) \geq 1$ then 
\begin{eqnarray}\label{ar2}
r_1(-\mathcal{M}(-\phi))&\geq& 0 \geq \Lambda
g(s,\mathcal{K}(s,\phi)) 
\geq   \Lambda r(\mathcal{K}
(s,\phi)-1)\\ &\geq&  \Lambda r(\exp
[\mathcal{M}(\phi)]-1) \geq r_1(\mathcal{M}(\phi)).\nonumber
\end{eqnarray}
Let now $b \in [0,0.5]$. In this case, the change of variables
$z(s) = - \ln x(t)$, where $s = s(t)$ reduces (\ref{gle}) to the form
\begin{equation}
\label{glema} z'(s) = - \Lambda g(s, \mathcal{K}(s,- z_s)),\quad \
z \in \R, \quad s\geq 0.
\end{equation}
We have, for $\mathcal{K} (s, -\phi) \geq 1$,
that 
\begin{eqnarray}\label{ar3}
r_2(\mathcal{M}(\phi)) &\leq& 0 \leq - \Lambda
g(s,\mathcal{K}(s,-\phi)) \leq
  -\Lambda
r(\mathcal{K} \{s, -\phi)-1\})\\ &\leq& 
 -\Lambda r(\exp
[\mathcal{M}(-\phi)]-1) \leq r_2(-\mathcal{M}(-\phi)),\nonumber
\end{eqnarray} where $r_2(x) = a\Lambda x/(1+(0.5-b)x)$. (Here we use
$w(x)=-r(e^{-x}-1)$ and argue as before). 
Next,
 if $\mathcal{K} (s, - \phi) \leq 1$,
we obtain 
\begin{eqnarray}\label{ar4}
r_2(-\mathcal{M}(- \phi)) &\geq& 0 \geq  - \Lambda
g(s,\mathcal{K}(s,-\phi(\cdot))) \geq 
-\Lambda r(\mathcal{K}(s, -\phi)-1)\\ &\geq& 
 -\Lambda r(\exp [-\mathcal{M}(\phi)]-1)
\geq r_2(\mathcal{M}(\phi)). \nonumber \end{eqnarray}

Hence, both equations (\ref{glem}) and (\ref{glema}) are of Carath\'eodory type and 
satisfy {\bf(GY)}. Now, if $\lim_{s\to +\infty}w(s) = w_*\not= 0$, then 
$\int_{\R_+}g(s,w(s))ds = \int_{\R_+}\lambda(t)f(t, w_1(t))dt$ for some $w_1(t)$ 
with $\lim_{t\to +\infty}w_1(t) = \exp(w_*) \not= 1$ so that {\bf(H3)} is also 
fulfilled.  Finally, the boundedness 
requirement from {\bf(H1)} follows easily from (\ref{ar1}), (\ref{ar2}),
(\ref{ar3}), (\ref{ar4}). In this way, we can finish the proof of Theorem
\ref{apl} by applying Theorem \ref{main} to Eqns. (\ref{glem}) and
(\ref{glema}).
\end{proof}

\begin{remark} \label{reu} Notice that Theorem \ref{apl} still remains true if, in its 
statement,  we replace the rational function $r(x)$ by any decreasing function 
$\tilde r: (-1,\infty) \to \R$ having  negative Schwarz derivative, 
satisfying (\ref{nf}), and such that $\tilde r'(0) = a < 0$ and $\tilde r$ is below bounded if
$-\tilde r''(0)/(2\tilde r'(0))\geq 1/2$. Indeed, in this case we still  can use Lemma \ref{l1}
while evaluating
$\rho(x) =
\tilde r(\exp(\pm x)-1))$ (see (\ref{ar1}), (\ref{ar2}),(\ref{ar3}), (\ref{ar4})). 
The sign of the second derivative of $\tilde r$ 
(associated with the sign of $b$ before) does not matter 
now since $\tilde r$ is defined for all $t > -1$.  
It is clear that $|a|\Lambda \leq 3/2$  implies the global stability 
 when $\rho''(0)\not = 0$, while for $\rho''(0)= 0$ 
we should assume that $|a|\Lambda < 3/2$. 
\end{remark}

\begin{example} Let us apply Theorem \ref{apl} to 
study the  food-limited population model \cite{gopalb, glq, kuang, liu, sy}
with variable, continuous and generally unbounded delay $h(t) \geq 0$ 
such that $\inf_{\R_+} (t-h(t)) = h_* $ is finite (obviously, $h_* \leq 0$): 
\begin{equation}
\label{mm}
N'(t) = \lambda (t)N(t) \frac{k - N^l(t-h(t))}{k + \nu(t)N^l (t-h(t))},
\ t\geq 0\ .
\end{equation}
Here $k, l > 0$, $\lambda \in \C([h_*, \infty), (0, \infty)),$
$\nu \in \C([0, \infty), [0, \infty))$  and $\int_{\R_+}\lambda = +\infty$. 

\begin{corollary}
\label{flm}
Let $\nu_0 = \inf_t \nu(t)\geq 0$, and
assume that for some $T\geq 0$ either $\nu_0\neq 1$ and $l\Lambda /(1+\nu_0) \le
\frac 32$, or
$\nu_0= 1$ and
$l\Lambda/(1+\nu_0) < \frac 32$, where $\Lambda =\sup_{t\geq
T}\left\{\int_{t-h(t)}^{t}
\lambda(s)ds\right\}$. 
Then every positive solution of equation
(\ref{mm}) converges to $k^{1/l}.$
\end{corollary}
\begin{proof}
With $N(t) =k^{1/l}\ x(s)$ and $s= s(t) = \Lambda
^{-1}\int_0^t\lambda(u)du$, Eq. (\ref{mm}) is  transformed 
into
$x'(s) = \Lambda x(s)f(s,
\mathcal{L}(s, x_s)),$
where $\mathcal{L}(s,\phi)=\phi(-\Lambda
^{-1}\int^{t(s)}_{t(s)-h(t(s))}\lambda(u)du)$, and $f(s,x)=(1-x^l)/(1+\nu(t(s))x^l)$.

In order to apply Theorem \ref{apl} to the {\it transformed} equation, 
we consider function
$$
f^*(s,x)=f(t(s),x+1)=\frac{1-(x+1)^l}{1+\nu(t(s))(x+1)^l}, \quad x > -1. 
$$

It is obvious that $xf^*(s,x)<0$ for all $x\neq 0$. Next, 
$(1-(x+1)^l)/(1+u(x+1)^l)$ is increasing with respect to $u$ for $x>0$ fixed
and decreasing with respect to $u$ for fixed $x<0$. Hence, $f^*(s,x)\geq \tilde r(x)$ for
$x>0$, and   $f^*(s,x)\leq \tilde r(x)$ for $x<0$, where
$$
\tilde r(x)=\frac{1-(x+1)^l}{1+\nu_0(x+1)^l}, \quad x > -1. 
$$

Finally, it is easy to check that  $\tilde r$ satisfies all 
conditions indicated in Remark \ref{reu}. Notice that $a=\tilde r'(0)=-l/(1+\nu_0)<0$
and  $\rho''(0) = 0$ only if $\nu_0=1$. 
\end{proof}

\begin{remark}
In the particular case when the delay is constant,  Corollary \ref{flm}
was proved in \cite{glq} under the stronger assumptions $l\Lambda/(1+\nu_0) \leq 1$ and 
\begin{equation}
\label{mmm}
\int_0^{\infty}\frac{\lambda(s)}{1+\nu(s)}ds=\infty. 
\end{equation}
Also, Liu in \cite{liu} considered 
\begin{equation}
\label{eqliu}
N'(t) = \lambda (t)N(t) 
\left(\frac{k - N(t-h)}{k + \nu(t)N(t-h)}
\right)^{\beta},
\ t \geq 0, 
\end{equation}
where $\beta = (2m+1)/(2n+1) \geq 1,\, m,n\in\N,$ getting  our $3/2$ condition only 
in the particular case when $\nu_0 =1$ (see  \cite[Theorem 2]{liu}). We remark that,
for
$\beta >1$,  the global attractivity of (\ref{eqliu}) can be always proved  once (\ref{mmm}) is
assumed. Notice also that the statement of Corollary 
\ref{flm} remains true if we replace all entries of  $N^l(t-h(t))$ in (\ref{mm}) 
by $N^l(\mu t)$, with $\mu \in (0,1)$, or by $ (\max_{u\in [t-h_0,t] } N(u))^l$ for 
some $h_0 > 0$ (thus we can include in our considerations equations 
with linearly transformed argument and equations with maxima). 
\end{remark}
\begin{remark} Assume (\ref{mmm}). In \cite{sy}, it has been 
established that the steady state of (\ref{mm}) with $h(t)= h$ 
is (locally) uniformly and asymptotically stable if
\begin{equation}
\label{mmmm}
l\int_{t-h}^{t}\frac{\lambda(s)}{1+\nu(s)}ds \leq \alpha < 3/2\, ,\, t\geq h. 
\end{equation}
This inequality is less restrictive than the one we obtained in Corollary \ref{flm}; 
thus, inspired by the Wright conjecture about the equivalence of global and local 
asymptotic stability, one can try to improve our result up to (\ref{mmmm}). 
However, it would be impossible: even with 
(\ref{mmmm}) satisfied, Eq. (\ref{mm}) can have nontrivial periodic solutions. 

\end{remark}

\end{example}


\begin{thebibliography}{30}

\bibitem{gopalb} K. Gopalsamy, ``Stability and oscillations in
delay differential equations of population
dynamics," Mathematics and its Applications, 74,
Kluwer, Dordrecht, 1992.

\bibitem{glq} E. A. Grove, G. Ladas and C. Qian, \textit{Global
attractivity in a `food limited' population model,} Dynamic Syst. Appl.,
\textbf{2} (1993), 243-250.

\bibitem{hl} J.~K.  Hale and  S.~M. Verduyn Lunel, ``Introduction to
functional differential equations", Applied Mathematical  Sciences,
Springer-Verlag, 1993.


\bibitem{ilt} A. Ivanov, E. Liz and S. Trofimchuk, \textit{Halanay inequality,
Yorke $3/2$ stability criterion, and differential equations with maxima,}
Tohoku Math. J. \textbf{54} (2002), 277-295.

\bibitem{kuang} Y. Kuang, ``Delay differential equations
with applications in population dynamics,"
Academic Press, 1993.

\bibitem{liu} Y. Liu, \textit{Global attractivity for a differential-difference
population model,} Appl. Math. E-Notes \textbf{1} (2001), 56--64 (electronic).

\bibitem{lprtt} E. Liz, M. Pinto, G. Robledo,
V. Tkachenko and S. Trofimchuk,  Wright type delay
differential equations with negative Schwarzian,
Discrete Contin. Dynam. Systems \textbf{9} (2003), 309--321.

\bibitem{ltt} E. Liz,
V. Tkachenko and S. Trofimchuk,  \textit{A global stability criterion for scalar
functional differential equations},
submitted (available at
{\tt http://arXiv.org/abs/math.DS/0112047).}

\bibitem{my} A.~D. Myshkis, ``Lineare Differentialgleichungen mit Nacheilendem
Argument," Deutscher Verlag Wiss., Berlin, 1955.

\bibitem{roy} H.~L. Royden, ``Real Analysis," The Macmillan Company, 1969.

\bibitem{sy} J.~W.-H. So and J.~S. Yu, \textit{
On the uniform stability for a `food-limited' population model with time delay,}
Proc. Roy. Soc. Edinburgh Sect. A  \textbf{125} (1995), 991--1002.


\bibitem{waltdcds} H.-O. Walther, \textit{Contracting return maps for  monotone
delayed feedback,} Discrete Contin. Dynam. Systems \textbf{7} (2001), 259--274.

\bibitem{wr} E.~M. Wright, \textit{A nonlinear difference-differential  equation,}
J. Reine Angew. Math. \textbf{194} (1955), 66--87.

\bibitem{yo}  J.~A. Yorke, \textit{Asymptotic stability for one dimensional
differential-delay
equations,}  J. Differential Equations \textbf{7}  (1970), 189--202.

\end{thebibliography}
\end{document}